\documentclass{amsart}

\usepackage{amssymb,amsthm,latexsym,amsfonts,amsmath}
\usepackage{url}
\usepackage{epsfig}
\usepackage[scanall]{psfrag}
\usepackage{subfigure}
\usepackage{algorithm2e}

\newtheorem{theorem}{Theorem}[section]

\newtheorem{proposition}[theorem]{Proposition}

\newtheorem{definition}[theorem]{Definition}

\DeclareMathOperator{\esssup}{ess\,sup}

\newcommand{\R}{{\mathbb{R}}}

\newcommand{\N}{{\mathbb{N}}}

\newcommand{\ie}{{\it i.e.}}
\newcommand{\eg}{{\it e.g. }}

\usepackage{color} % Tabuada

\begin{document}

\title[Decentralized Event-triggered Control over WSAN]{Decentralized event-triggered control\\ over wireless sensor/actuator networks}
\author{Manuel Mazo Jr and Paulo Tabuada}
\thanks{M. Mazo Jr is with INCAS$^3$, Assen and the Department of Discrete Technology and Production Automation, University of Groningen, The Netherlands, 
{\tt\small m.mazo@rug.nl}\\
P. Tabuada is with the Department of Electrical Engineering, University of California, Los Angeles, CA 90095-1594,{\tt\small tabuada@ee.ucla.edu}%
}
%\author[Manuel Mazo Jr.]{Manuel Mazo Jr.$^1$} 
%\author[Paulo Tabuada]{Paulo Tabuada$^1$}
%\address{$^1$Department of Electrical Engineering\\
%University of California at Los Angeles,
%Los Angeles, CA 90095}
%\email{\{mmazo, tabuada\}@ee.ucla.edu}
%\urladdr{http://www.ee.ucla.edu/~mmazo}
%\urladdr{http://www.ee.ucla.edu/~tabuada}

\begin{abstract}
In recent years we have witnessed a move of the major industrial automation providers into the wireless domain. While most of these companies already offer wireless products for measurement and monitoring purposes, the ultimate goal is to be able to close feedback loops over wireless networks interconnecting sensors, computation devices, and actuators. 
In this paper we present a decentralized event-triggered implementation, over sensor/actuator networks, of centralized nonlinear 
controllers. Event-triggered control has been recently proposed as an alternative to the more traditional 
periodic execution of control tasks. In a typical event-triggered implementation, the control signals are kept constant until the violation of a condition on the state of the plant triggers the re-computation of the control signals. 
The possibility of reducing the number of re-computations, and thus of transmissions, while guaranteeing desired levels of control performance, makes event-triggered control very appealing in the context of sensor/actuator networks. 
In these systems the communication network is a shared resource and event-triggered
implementations of control laws offer a flexible way to reduce network utilization. Moreover reducing the number of times that a feedback control law is executed implies a reduction in transmissions and thus a reduction in energy expenditures of battery powered wireless sensor nodes.
\end{abstract}

\maketitle

\section{Introduction}

%A large amount of the research available on sensor networks has focused on the acquisition of measurements from physical processes. Many of the developed applications concentrate on how to obtain this information for posterior off-line analysis~\cite{App1,App2}. Other researchers are concerned with on-line processing of this information for different applications such as tracking~\cite{OhSastry,SinopoliSastry}, distributed optimization~\cite{RabbatNowak} or mapping~\cite{Map}. 
%All of these applications show a common desire for small power consumption in order to extend 
%the life span of the system.  We address the problem of reducing power consumption by reducing communication in applications in which actuation plays a major role, namely control applications. 

%The use of wireless sensor/actuator networks (WSAN) for control purposes has been receiving attention in academia for quite some time. Industry has recently showed interest as evidenced by the WirelessHART initiative~\cite{wHART_URL}.

For many years, control engineers have designed their controllers as if there were infinite-bandwidth, 
noise- and delay-free channels between sensors, controllers, and actuators. The effects of 
non-idealities in the channels, in practice, could be mitigated by employing better hardware. 
However, on implementations over Wireless Sensor Actuator Networks (WSAN) these limitations of the communication medium can no longer be neglected. This fact, combined with the recent interest from industry, \eg the WirelessHART initiative~\cite{wHART_URL}, have fueled the study of control under communication constraints in the past decade. Much research has been devoted to the effects of: quantization in the sensors; delay and jitter; limited bandwidth; or even packet losses. Some good overviews of these topics can be found in the report resulting from the RUNES project~\cite{RUNES}, and the special issue of the IEEE proceedings~\cite{special_proceedings}.

%A different approach to save communication energy in control over WSAN's was presented by Rozell  and Johnson in~\cite{Rozell}, but still assuming periodic communications. Based on the redundancy of the information they schedule the amount of power and bits each node is assigned, resulting in a more efficient use of the available energy.
One aspect common to most modern control systems, and something assumed in most of the studies 
mentioned above, is the implementation of control strategies in embedded microprocessors.
But in controlling the physical world, which is of continuous nature, the use of 
microprocessors brings a new question: how often should we sample the physical environment~\cite{Franklin}? 
Many researchers have worked on the analysis of this sole problem. Tools like
the delta-transform~\cite{Middleton} were developed, and many
books discussed this issue~\cite{Goodwin,Houpis}. More recently, Nesic and collaborators have proposed techniques to select periods retaining closed-loop stability in networked systems~\cite{Nesic:2001p2652, Nesic:2009}. However, engineers still rely mostly on rules of
thumb such as sampling with a frequency 20
times the system bandwidth, and then check if it actually works~\cite{Franklin,Goodwin,Houpis}.
A shift in perspective was brought by the notion of event-triggered control~\cite{arzen99}, \cite{astrom}. In event-triggered control, instead of periodically updating the control input, the update instants are generated by the violation of a condition on the state of the plant. 
Many researchers have shown a renewed interest on these techniques~\cite{Heemels08,Cervin08,Rabi08,Rabi082,Molin10,Lunze10}.
Recently, one of the authors proposed a formalism to generate asymptotically stable event-triggered implementations of nonlinear controllers~\cite{Tabuada}, and in~\cite{MazoJr:2008p1340} the authors explored the application of event-triggered and self-triggered techniques to distributed implementations of linear controllers. For more details about these event-triggered and self-triggered techniques we refer the reader to~\cite{Anta:2008p2867} and~\cite{MazoAnta}.
Following the formalism in~\cite{Tabuada}, Wang and Lemmon proposed a distributed event-triggered implementation for weakly-coupled  distributed systems~\cite{WangLemmon09}.
The present work complements the techniques described in~\cite{WangLemmon09} by addressing systems without weak-coupling assumptions.

The main contribution of this paper is a strategy for the construction of decentralized event-triggered implementations over WSAN of centralized controllers. The event-triggered techniques introduced in~\cite{Tabuada}
are based on a criterion that depends on the norm of the vector of
measured quantities. This is natural in the setting discussed in~\cite{Tabuada}
since sensors were collocated with the micro-controller. However, in a WSAN the physically distributed sensor nodes do not have access to all the measured quantities. Hence, we
cannot use the same criterion to determine when the control signal should be re-computed. Using classical observers or estimators (as the Kalman filter) would require
filters of dimension as large as the number of states in each sensor node, which would be  unpractical given the low computing capabilities of sensor nodes. Moreover, we do not assume observability from every measured output, thus ruling out observer-based techniques. Approaches based on consensus algorithms are also unpractical as they require large amounts of communication and thus large energy expenditures by the sensor nodes. Instead, we present an approach to decentralize a centralized event-triggered condition that relies only on the locally measured quantities. Our technique also provides a mechanism to enlarge the resulting times between controller re-computations without altering performance guarantees.

We do not address in this paper practical issues such as delays or jitter in the communication and focus solely on the reduction of the actuation frequency (with its associated communication and energy savings). In particular, the issue of communication delays has been shown to be easily addressed in the context of event-triggered control in~\cite{Tabuada} and similarly in~\cite{WangLemmon09}. The approach followed in those papers is applicable to the techniques introduced in this paper. Moreover, our techniques can be implemented over the WirelessHART standard~\cite{wHART_URL}, which addresses other communication concerns such as medium access control, power control, and routing. 

The present paper is organized as follows: we introduce basic notation in Section~\ref{sec:notation};
Section~\ref{sec:decentralized_ET} states the problem, briefly reviews the results of~\cite{Tabuada} and presents our proposal for decentralization; the paper finalizes with an example in Section~\ref{sec:example} and a discussion in Section~\ref{sec:discussion}.

\section{Notation}
\label{sec:notation}

We denote by $\N$ the natural numbers, by \mbox{$\N_0=\N\cup\lbrace0\rbrace$}, 
by $\R^+$ the positive real numbers,  and
by $\R_0^+=\R^+\cup\lbrace0\rbrace$. The usual Euclidean ($l_2$) vector norm is
represented by $|\cdot|$. When applied to a matrix, $|\cdot|$ denotes the
$l_2$ induced matrix norm. A matrix
$P\in\R^{n\times n}$ is said to be positive definite, denoted $P>0$, whenever
$x^TPx>0$ for all $x\neq 0$, $x\in\R^n$. By $\lambda_m(P),\lambda_M(P)$
we denote the minimum and maximum eigenvalues of $P$ respectively. A function
$\gamma:[0,\infty[\to\R^+_0$, is of class~$\mathcal{K}_\infty$ if it is
continuous, strictly increasing, $\gamma(0)=0$ and $\gamma(s)\to\infty$ as
$s\to\infty$. 
Given an essentially bounded function $\xi:\R^+_0\to\R^n$ we denote by $\Vert
\xi \Vert_\infty$ its $\mathcal{L}_\infty$ norm, \ie, $\Vert \delta
\Vert_\infty=\esssup_{t\in\R^+_0}\lbrace{|\xi(t)|\rbrace}<\infty$.

In the following we consider systems defined by differential equations of
the form:
\begin{equation}
\label{eq:control_sys}
\frac{d}{dt}\xi=f(\xi,\upsilon)
\end{equation}
with input $\upsilon:\R^+_0\to\R^m$ an essentially bounded piecewise continuous function of
time and $f:\R^n\times\R^m\to\R^n$ a smooth map. We also use the simpler
notation $\dot{\xi}=f(\xi,\upsilon)$ to refer to~(\ref{eq:control_sys}).
We refer to such systems as {\it control systems}. Solutions
of~(\ref{eq:control_sys}) with initial condition $x$ and input $\upsilon$,
denoted by $\xi_{x\upsilon}$, satisfy: $\xi_{x\upsilon}(0)=x$ and
$\frac{d}{dt}\xi_{x\upsilon}(t)=f(\xi_{x\upsilon}(t),\upsilon(t))$
for almost all $t\in\R^+_0$. The notation will be relaxed by dropping
the subindex when it does not contribute to the clarity of exposition. 
A feedback law for a control system is a smooth map
$k:\R^n\to\R^m$; we sometimes refer to such a law as a {\it controller}
for the system.

\section{Decentralized event-triggered control}
\label{sec:decentralized_ET}

Consider a nonlinear control system and a hardware platform consisting of a set of wireless sensors and actuators and a computation node. This last node is in charge of computing the control signal with the measurements obtained from the sensors. We consider scenarios in which none of these sensor nodes has access to the full state of the plant.
We model the execution of the control loop in three steps: data retrieval from sensors, controller computation, and provision of the control commands to the actuators.
Furthermore, we assume that the computation of the controller happens in just one device which retrieves all the measurement information from the sensors, computes the inputs for all actuators, and disseminates these new commands to the actuator nodes. This scenario is a typical configuration considered in the WirelessHART standard, see~\cite{Soldati}, which addresses the problem of scheduling links and channels for disseminating the information in WirelessHART networks.

Our goal is to provide a mechanism triggering the execution of the control loop which reduces the frequency of the controller updates. In order to reduce the frequency of controller updates we abandon the periodic transmission paradigm, and instead we propose to close the loop whenever certain events happen. In particular, we consider the event-triggered implementation techniques proposed in~\cite{Tabuada} which guarantee the asymptotic stability of the closed-loop system. These techniques, however, require the knowledge of the full state to decide when to trigger new updates, but such information is not available at any sensing node under our premises. In the following we discuss a decentralization of the decision process triggering controller updates. We propose the use of conditions depending solely on the information available at each node. Whenever any of these conditions is violated at a node, this node informs the computation device. Upon receipt of such an event, the computation device requests fresh measurements, updates the control signals, and forwards the new commands to the actuation nodes.

\subsection{Event-triggered control}

We begin by revisiting the results from~\cite{Tabuada}, which serve as the basis for the rest of our work.  Let us start by considering a nonlinear control system:
\begin{equation}
\label{eq:dynamics}
\dot{\xi}=f(\xi,\upsilon)
\end{equation}
and assume that a feedback control law $k:\R^n\to\R^m$, 
\mbox{$\upsilon=k(\xi)$} is available, rendering the closed-loop system:
\begin{equation}
\label{eq:closed_loop_cont}
\dot{\xi}=f(\xi,k(\xi+\varepsilon))
\end{equation}
\emph{input-to-state stable} (ISS)~\cite{Sontag:2008p2453} with respect to measurement errors $\varepsilon:\R_0^+\to\R^n$.  We do not provide the definition of ISS, but rather the following characterization that lies at the heart of our techniques:
\begin{definition}
A smooth function $V:\R^n\to\R_0^+$ is said to be an ISS Lyapunov function for the closed-loop system~(\ref{eq:closed_loop_cont}) if there exists class $\mathcal{K}_\infty$ functions $\underline{\alpha}$,$\overline{\alpha}$, $\alpha$ and $\gamma$ such that for all $x\in\R^n$ and $e\in\R^n$ the following is satisfied:
\begin{eqnarray}
\label{eq:V_bounded}\notag
\underline{\alpha}(|x|)\leq& V(x)&\leq \overline{\alpha}(|x|)\\
\label{eq:Vdot_bounds}
\frac{\partial V}{\partial x} f(x,k(x+e))&\leq& -\alpha(|x|)+\gamma(|e|).
\end{eqnarray}
\end{definition}
The closed-loop system~(\ref{eq:closed_loop_cont}) is said to be ISS with respect to measurement errors $\varepsilon$, if there exists an ISS Lyapunov function for~(\ref{eq:closed_loop_cont}). 

In a sample-and-hold implementation of the control law $k(\xi)$, the input signal is held constant between update times, \ie:
\begin{eqnarray}\notag
\label{eq:HybSys1}
\dot{\xi}(t)&=&f(\xi(t),\upsilon(t))\\
\upsilon(t)&=&k(\xi(t_k)), \; t\in[t_k,t_{k+1}[,
\end{eqnarray}
where $\lbrace{t_k\rbrace}_{k\in\N^+_0}$ is a divergent sequence of update
times. An event-triggered implementation defines such a
sequence of update times $\lbrace{t_k\rbrace}_{k\in\N^+_0}$ for the controller,
rendering the closed loop system asymptotically stable.

We now consider the signal $\varepsilon:\R_0^+\to\R^n$ defined by \mbox{$\varepsilon(t)=\xi(t_k)-\xi(t)$} for $t\in[t_k,t_{k+1}[$ and regard it as a measurement error. By doing so, we can rewrite~(\ref{eq:HybSys1}) 
for $t\in[t_k,t_{k+1}[$ as:
\begin{eqnarray}\notag
\label{eq:HybSys2}
\dot{\xi}(t)&=&f(\xi(t),k(\xi(t)+\varepsilon(t))),\\ \notag
\dot{\varepsilon}(t)&=&-f(\xi(t),k(\xi(t)+\varepsilon(t))),\;\varepsilon(t_k)=0.
\end{eqnarray}
Hence, as~(\ref{eq:closed_loop_cont}) is ISS with respect to measurement errors $\varepsilon$, from~(\ref{eq:V_bounded}) we know that by enforcing:
\begin{equation}
\label{eq:ineq_orig}
\gamma(|\varepsilon(t)|)\leq\rho\alpha(|\xi(t)|),\;\forall t>0,\; \rho\in]0,1[
\end{equation}
the following holds:
\begin{equation}\notag
\frac{\partial V}{\partial x} f(x,k(x+e))\leq -(1-\rho)\alpha(|x|),\;\forall\,x,e\in\R^n
\end{equation}
and asymptotic stability of the closed-loop follows. 
Moreover, if one assumes that the system operates in some compact set  $S\subseteq \R^n$ and $\alpha^{-1}$ and $\gamma$ are Lipschitz continuous on $S$, the inequality~(\ref{eq:ineq_orig}) can be replaced by the simpler inequality $|\varepsilon(t)|^2\leq\sigma|\xi(t)|^2$, for a suitably chosen $\sigma>0$. 
Hence, if the sequence of update times $\lbrace{t_k\rbrace}_{k\in\N^+_0}$ is such that:
\begin{equation}
\label{eq:condition}
|\varepsilon(t)|^2\leq\sigma |\xi(t)|^2,\quad t\in [t_k,t_{k+1}[,
\end{equation}
the sample-and-hold implementation~(\ref{eq:HybSys1}) is guaranteed to render the closed loop system asymptotically stable. 

Condition~(\ref{eq:condition}) defines an event-triggered implementation that consists of continuously checking~(\ref{eq:condition}) and triggering the recomputation of the control law as soon as the inequality evaluates to equality. Note that recomputing the controller at time $t=t_k$ requires a new state measurement and thus resets the error $\varepsilon(t_k)=\xi(t_k)-\xi(t_k)$ to zero which enforces~(\ref{eq:condition}).

\subsection{Decentralized event-triggering conditions}

We consider, for simplicity of presentation, a decentralized scenario in which each state variable is measured by a different sensor. However, the same ideas apply to more general decentralized scenarios as we briefly discuss at the end of Section~\ref{ssec:adapt}.
In this setting, no sensor can evaluate condition~(\ref{eq:condition}), since~(\ref{eq:condition}) requires the knowledge of the full state vector $\xi(t)$.
Our goal is to provide a set of simple conditions that each sensor can check locally to decide when to trigger a controller update, thus triggering also the transmission of fresh measurements from sensors to the controller. 

Using a set of parameters $\theta_1, \theta_2,\ldots, \theta_n\in\R$ such that $\sum_{i=1}^{n} \theta_i=0$, we can rewrite inequality~(\ref{eq:condition}) as:
$$
\sum_{i=1}^n\left( \varepsilon_i^2(t)-\sigma \xi_i^2(t)\right)\leq 0 = \sum_{i=1}^n \theta_i,
$$
where $\varepsilon_i$ and $\xi_i$ denote the $i$-th coordinates of $\varepsilon$ and $\xi$ respectively.
Hence, the following implication holds:
\begin{equation}
\label{eq:imply_cond}
\bigwedge_{i=1}^n  \; \left( \varepsilon_i^2(t)-\sigma \xi_i^2(t)\leq \theta_i \right) \Rightarrow |\varepsilon(t)|^2\leq\sigma |\xi(t)|^2,
\end{equation}
which suggests the use of:
\begin{equation}
\label{eq:decentralized_conditions}
\varepsilon_i^2(t)-\sigma \xi_i^2(t)\leq \theta_i
\end{equation} 
as the local event-triggering conditions. 

In this decentralized scheme, whenever any of the local conditions~(\ref{eq:decentralized_conditions}) becomes an equality, the controller is recomputed. We denote by $t_k+\tau_i(x)$ the first time at which~(\ref{eq:decentralized_conditions}) is violated, when $\xi(t_k)=x$, $\varepsilon(t_k)=0$. If the time elapsed between two events triggering controller updates is smaller than the minimum time $\tau_{min}$ between updates of the centralized event-triggered implementation\footnote{It was proved in~\cite{Tabuada} that such a minimum time exists for the centralized condition, and that lower bounds can be explicitly computed.}, the second event is ignored and the controller update is scheduled $\tau_{min}$ units of time after the previous update. 

Not having an equivalence in~(\ref{eq:imply_cond}) entails that this decentralization approach is in general conservative: times between updates will be shorter than in the centralized case. The vector of parameters $\theta=[\theta_1\,\theta_2\,\ldots\theta_n]^T$ can be used to reduce the mentioned conservatism and thus reduce utilization of the communication network.
It is important to note that the vector $\theta$ can change every time the control input is updated. From here on we show explicitly this time dependence of $\theta$ by writing $\theta(k)$ to denote its value between the update instants $t_k$ and $t_{k+1}$. 
Following the presented approach, as long as $\theta$ satisfies $\sum_{i=1}^{n} \theta_i(k)=0$, the stability of the closed-loop is guaranteed regardless of the specific value that $\theta$ takes and the rules used to update $\theta$.

We summarize the previous discussion in the following proposition:
\begin{proposition}
For any choice of $\theta$ satisfying:
$$\sum_{i=1}^{n} \theta_i(k)=0,\;\forall\, k\in\N^+_0,$$
the sequence of update times $\lbrace{t_k\rbrace}_{k\in\N^+_0}$ given by:
%\begin{small} 
\setlength{\arraycolsep}{0.3pt}
\begin{eqnarray}\notag
t_{k+1}&=&t_k+\max \lbrace \tau_{min}, \min_{i=1,\ldots,n} \tau_i(\xi(t_k))\rbrace\\ \notag
\tau_i(\xi(t_k))&=&\min\lbrace \tau\in\R_0^+\, |\; \epsilon_i^2(t_k+\tau)-\sigma\xi_i^2(t_k+\tau)=\theta_i(k) \rbrace
\end{eqnarray}
%\end{small}
renders the system~(\ref{eq:HybSys1}) asymptotically stable.
%$$\varepsilon_i^2(t)-\sigma \xi_i^2(t)\leq \theta_i(k),\; \forall\, i\in\lbrace{1, 2, \ldots,n\rbrace}, \,t\in[t_k,t_{k+1}[,\; k\in\N^+_0,$$
\end{proposition}

\subsection{Decentralized event-triggering with on-line adaptation}
\label{ssec:adapt}

We present now a family of heuristics to adjust the vector $\theta$ whenever the control input is updated. 
We define the \emph{decision gap} at sensor $i$ at time $t\in[t_k,t_{k+1}[$ as: 
$$G_i(t)= \varepsilon_i^2(t)-\sigma \xi_i^2(t)- \theta_i(k).$$ 
The heuristic aims at equalizing the decision gap at some future time. 
We propose a family of heuristics parametrized by an \emph{equalization time} $t_e$ and an \emph{approximation order} $q$.

For the \emph{equalization time} $t_e:\N_0\to\R^+$ we present the following two choices:
constant and equal to the minimum time between controller updates $t_e(k)=\tau_{min}$;
the previous time between updates $t_e(k)=t_k-t_{k-1}$.

The \emph{approximation order} is  the order of the Taylor expansion used to estimate the decision gap at the equalization time $t_e$:
$$
\hat{G}_i(t_k+t_e)=\hat{\varepsilon}_i^2(t_k+t_e)-\sigma \hat{\xi}_i^2(t_k+t_e)- \theta_i(k).
$$
where for $t\in[t_{k},t_{k+1}[$:
\begin{small}
\begin{eqnarray*}
\hat{\xi}_i(t)&=&\xi_i(t_k)+\dot{\xi}_i(t_k)(t-t_k)+\frac{1}{2}\ddot{\xi}_i(t_k)(t-t_k)^2+\ldots\\&&+\frac{1}{q!}\xi^{(q)}_i(t_k)(t-t_k)^q,\\
\hat{\varepsilon}_i(t)&=&\phantom{x_i}0\phantom{t_k}-\dot{\xi}_i(t_k)(t-t_k)-\frac{1}{2}\ddot{\xi}_i(t_k)(t-t_k)^2-\ldots\\&&-\frac{1}{q!}\xi^{(q)}_i(t_k)(t-t_k)^q,
\end{eqnarray*}
\end{small}
using the fact that $\dot{\varepsilon}=-\dot{\xi}$ and $\varepsilon(t_k)=0$. 

Finally, once an equalization time $t_e$ and an approximation order $q$ are chosen, the vector $\theta(k)\in\R^n$ is computed so as to satisfy:
\begin{small}
\begin{eqnarray}\notag
&&\hat{G}_i(t_k+t_e)=\hat{G}_j(t_k+t_e)\qquad\forall i,j \in \lbrace1,2,\ldots,n\rbrace,\\\label{eq:adaptation}
&&\sum_{i=1}^n\theta_i(k) =0.
\end{eqnarray}
\end{small}
Note that finding such $\theta$, after the estimates $\hat{\xi}$ and $\hat{\varepsilon}$ have been computed, amounts to solving a system of $n$ linear equations.
% \begin{scriptsize}
% \begin{eqnarray}
% \label{eq:adaptation}
% &\left[
% \begin{array}{cccccc}
% 1 & -1 & 0 & 0 & \ldots & 0\\
% 0 &  1 & -1 & 0 & \ldots & 0\\
% 0 & 0 & \ddots & \ddots & 0 & 0\\
% 0 & 0 & 0 & \ldots & 1 & -1\\
% 1 & 1 & 1 & \ldots & 1 & 1
% \end{array}
% \right]
% \left[
% \begin{array}{c}
% \theta_1(k)\\
% \theta_2(k)\\
% \vdots\\
% \theta_{n-1}(k)\\
% \theta_n(k)
% \end{array}
% \right]=
% \left[
% \begin{array}{c}
% \delta_{12}(t_k+t_e)\\
% \delta_{23}(t_k+t_e)\\
% \vdots\\
% \delta_{(n-1)n}(t_k+t_e)\\
% 0
% \end{array}
% \right],\\\notag
% \\\notag
% &\delta_{ij}(t)=\left(\hat{\varepsilon}_i^2(t)-\sigma \hat{\xi}_i^2(t)\right)-\left(\hat{\varepsilon}_j^2(t)-\sigma \hat{\xi}_j^2(t)\right).
% \end{eqnarray}
% \end{scriptsize}
Note also that $\theta$ is computed\footnote{The resulting $\theta$ computed in this way could be such that for some sensor $i$, $-\xi_i^2(t_k)>\theta_i(k)$. Such choice of $\theta$ results in an immediate violation of the triggering condition at $t=t_k$, \ie, $\tau_i(\xi(t_k))$ would be zero. In practice, when the unique solution of~(\ref{eq:adaptation}) results in $-\xi_i^2(t_k)>\theta_i(k)$, one resets $\theta$ to some default value such as the zero vector.} in the controller node, which has access to $\xi(t_k)$.

The choice of $t_e$ and $q$ has a great impact on the amount of actuation required. The use of a large $t_e$ leads, in general, to poor estimates of the state of the plant at time $t_k+t_e$ and thus degrades the equalization of the gaps. On the other hand, one expects that equalizing at times $t_k+t_e$ as close as possible to the next update time $t_{k+1}$ (according to the centralized event-triggered implementation) provides larger times between updates. In practice, these two objectives (small $t_e$, and $t_{k}+t_e$ close to the ideal $t_{k+1}$) can be contradictory, namely when the time between controller updates is large. The effect of the order of approximation $q$ depends heavily on $t_e$ and enlarging $q$ does not necessarily improve the estimates. An heuristic providing good results in several case studies performed by the authors is given by Algorithm~\ref{alg:heuristic}.
\begin{small}
\begin{algorithm}
%\SetAlgoLined
\SetLine
\KwIn{$q$, $t_{k-1}$, $t_k$, $\tau_{min}$, $\xi(t_k)$}
\KwOut{$\theta(k)$}
$t_e:=t_k-t_{k-1}$\;
Compute $\theta(k)$ according to equation~(\ref{eq:adaptation})\;
\If{$\exists\, i\in \lbrace1,2,\ldots,n\rbrace$ such that $-\xi_i^2(t_k)>\theta_i(k)$}{
$t_e:=\tau_{min}$\;
Compute $\theta(k)$ according to equation~(\ref{eq:adaptation})\;
\If{$\exists\, i\in \lbrace1,2,\ldots,n\rbrace$ such that $-\xi_i^2(t_k)>\theta_i(k)$}{
$\theta(k):=0$\;
}
}
\caption{The $\theta$-adaptation heuristic algorithm.\label{alg:heuristic}}
\end{algorithm}
\end{small}

While we assumed, for simplicity of presentation, that each node measured a single state of the system, in practice there may be scenarios in which one sensor has access to several (but not all) states of the plant. The same approach applies by considering local triggering rules of the kind $|\bar\varepsilon_i(t)|^2 -\sigma|\bar\xi_i(t)|^2\leq \theta_i$, where $\bar\xi_i(t)$ is now the vector of states sensed at node $i$, $\bar\varepsilon_i(t)$ is its corresponding error vector, and $\theta_i$ is a scalar.
%Note, however, that to guarantee the existence of a minimum time between controller updates as in~\cite{Tabuada}, one needs to update all actuators at the same time. Nonetheless, this requirement could be relaxed in particular cases, under suitable conditions for the Lyapunov function.

\subsection{Comments on practical implementations}

The proposed technique, while clearly reducing the amount of information that needs to be transmitted from sensors to actuators, might suggest that sensor nodes need to be continuously listening for events triggered at other nodes. This poses a practical problem since the energy required to keep the radios on to listen for possible events could potentially be very large. In practice, however, sensor nodes have their radio modules asleep most of the time and are periodically awaken according to a time multiplexing medium access protocol. Time multiplexing is typically used in protocols for control over wireless networks, like WirelessHART, in order to provide non-interference and strict delay guarantees. The use of time multiplexing can be accommodated in the proposed technique by regarding its effect as a bounded and known delay between the generation of an event and the corresponding change in the control signal. As was shown in~\cite{Tabuada}, delays can be accommodated in event-triggered implementations by adequately reducing the value of $\sigma$, therefore making the triggering conditions more conservative.

\section{Examples and simulation results}
\label{sec:example}

\begin{figure}[ht]
\begin{center}
\includegraphics[scale=0.92]{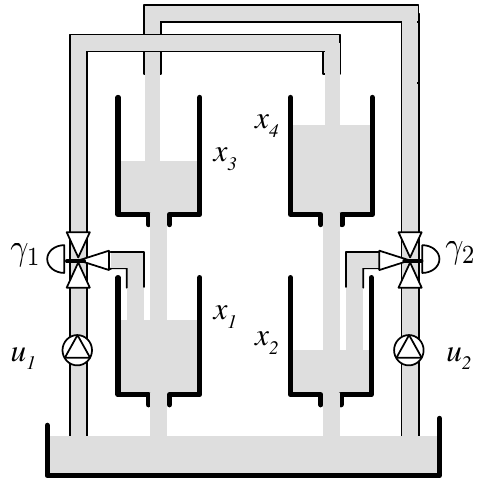}
\caption{The quadruple-tank system.}
\label{fig:tanks}
\end{center}
\end{figure}

We present in what follows an example illustrating the effectiveness of the proposed technique.
We select the quadruple-tank model from~\cite{allgoewer} describing the multi-input multi-output nonlinear system consisting of four water tanks as shown in Figure~\ref{fig:tanks}. The water flows from tanks $3$ and $4$ into tanks $1$ and $2$, respectively, and from these two tanks to a reservoir. 
The state of the plant is composed of the water levels of the tanks: $\xi_1$, $\xi_2$, $\xi_3$ and $\xi_4$. Two inputs are available: $\upsilon_1$ and $\upsilon_2$, the input flows to the tanks. The input flows are split at two valves $\gamma_1$ and $\gamma_2$ into the four tanks. The positions of these valves are given as parameters of the plant. The goal is to stabilize the levels $x_1$ and $x_2$ of the lower tanks at some specified values $x^*_1$ and $x^*_2$.

The system dynamics are given by the equation:
\begin{equation}\notag
\dot\xi(t)=f(\xi(t))+g_c \upsilon,
\end{equation}
with:
\begin{small}
\begin{equation}\notag
f(x)=\left[
\begin{array}{c}
-\frac{a_1\sqrt{2gx_1}}{A1}+\frac{a_3\sqrt{2gx_3}}{A1}\\
-\frac{a_2\sqrt{2gx_2}}{A2}+\frac{a_4\sqrt{2gx_4}}{A2}\\
-\frac{a_3\sqrt{2gx_3}}{A3}\\
-\frac{a_4\sqrt{2gx_4}}{A4}
\end{array}
\right],\;
g_c=\left[
\begin{array}{cc}
\frac{\gamma_1}{A1} & 0\\
0 & \frac{\gamma_2}{A2}\\
0 & \frac{1-\gamma_2}{A3}\\
\frac{1-\gamma_1}{A4} & 0
\end{array}
\right],
\end{equation}
\end{small}
and $g$ denoting gravity's acceleration and $A_i$ and $a_i$ denoting the cross sections of the $i-th$ tank and outlet hole respectively.

The controller design from~\cite{allgoewer} requires the extension of the plant with two extra artificial states $\xi_5$ and $\xi_6$. These states are nonlinear integrators used by the controller to achieve zero steady-state offset and evolve according to:
\begin{eqnarray}\notag
\dot\xi_5(t)=k_{I1}a_1\sqrt{2g}\left(\sqrt{\xi_1(t)}-\sqrt{x^*_1}\right),\\\notag
\dot\xi_6(t)=k_{I2}a_2\sqrt{2g}\left(\sqrt{\xi_2(t)}-\sqrt{x^*_2}\right),
\end{eqnarray}
where $k_{I1}$ and $k_{I2}$ are design parameters of the controller.
Note how stabilizing the extended system implies that in steady-state $\xi_1$ and $\xi_2$ converge to the desired values $x^*_1$ and $x^*_2$. We assume in our implementation that the sensors measuring $\xi_1$ and $\xi_2$, also compute $\xi_5$ and $\xi_6$ locally and sufficiently fast. Hence, we can consider $\xi_5$ and $\xi_6$ as regular state variables.

The controller proposed in~\cite{allgoewer} is given by the following feedback law:
\begin{equation}
\upsilon(t)=-K(\xi(t)-x^*)+u^*
\end{equation}
with
%\begin{footnotesize}
\begin{small}
\begin{eqnarray}
\notag \label{eq:u_star}
u^*&=&\left[
\begin{array}{cc}
\gamma_1 & 1-\gamma_2\\
1-\gamma_1 & \gamma_2
\end{array}
\right]^{-1}\left[\begin{array}{c}
a_1\sqrt{2gx^*_1}\\
a_2\sqrt{2gx^*_2}
\end{array}
\right]\\&=&\left[
\begin{array}{cc}
0 & 1-\gamma_2\\
1-\gamma_1 & 0
\end{array}
\right]^{-1}\left[\begin{array}{c}
a_1\sqrt{2gx^*_3}\\
a_2\sqrt{2gx^*_4}
\end{array}
\right],
\end{eqnarray}
\end{small}
and $K=QP$ where $Q$ is a positive definite matrix and $P$ is given by
%\begin{footnotesize}
\setlength{\arraycolsep}{0.3pt}
\begin{scriptsize}
\begin{equation}\notag
%P=\left[
%\begin{array}{cc}
%\gamma_1k_1 &  (1-\gamma_2)k_1 \\
%(1-\gamma_1)k_2 & \gamma_2k_2 \\
%0 &  (1-\gamma_2)k_3 \\
%(1-\gamma_1)k_4 &  0 \\
%\gamma_1k_1 &  (1-\gamma_2)k_1 \\
%(1-\gamma_1)k_2 &  \gamma_2k_2
%\end{array}
%\right]^T.
P=\left[
\begin{array}{cccccc}
\gamma_1k_1 &(1-\gamma_1)k_2 &0 & (1-\gamma_1)k_4 &\gamma_1k_1 &(1-\gamma_1)k_2\\
(1-\gamma_2)k_1 &\gamma_2k_2 &(1-\gamma_2)k_3 &0 &(1-\gamma_2)k_1&\gamma_2k_2
\end{array}
\right],
\end{equation}
\end{scriptsize}
%\end{footnotesize}
where $k_1$, $k_2$, $k_3$ and $k_4$ are design parameters of the controller.
Note how equation~(\ref{eq:u_star}) can be used to compute $x^*_3$ and $x^*_4$ from the specified $x^*_1$ and $x^*_2$. 
When computing the control $\upsilon$, the remaining entries $x^*_5$ and $x^*_6$ of $x^*=[x^*_1\;x^*_2\;x^*_3\;x^*_4\;x^*_5\;x^*_6]^T$  can be set to any arbitrary (fixed) values $\hat{x}^*_5$ and $\hat{x}^*_6$. This can be done because the errors: $\hat{x}^*_5-x^*_5$ and $\hat{x}^*_6-x^*_6$, between the arbitrary values and the actual states $x_5^*$ and $x_6^*$ of the equilibrium, can be reinterpreted as a perturbation on the initial states $\xi_5(0)$ and $\xi_6(0)$.

%The provided controller renders the closed loop in the form:
%\begin{equation}\notag
%\dot{\xi}(t)=(J-R)\nabla H_d(\xi(t))
%\end{equation}
%with
%\begin{equation}
%J-R=\left[
%\begin{array}{cccccc}
%\frac{-1}{k_1A_1} &  0 & \frac{1}{k_3A_1} & 0 & 0 & 0\\
%0 & \frac{-1}{k_2A_2} & 0 & \frac{1}{k_4A_2} & 0 & 0\\
%0 &  0 & \frac{-1}{k_3A_3} & 0 & 0 & 0\\
%0 &  0 & 0 & \frac{-1}{k_4A_4} & 0 & 0\\
%\frac{k_{I1}}{k_1} &  0 & 0 & 0 & \frac{-k_{I1}}{k_1} & 0\\
%0 &  \frac{k_{I2}}{k_2} & 0 & 0 & 0 & \frac{-k_{I2}}{k_2}
%\end{array}
%\right],
%\end{equation}
%\begin{eqnarray}\notag
%H_d(x)&=&\frac{1}{2}(x-x^*)^TP^TQP(x-x^*)-u^{*T}Px+\sum_{i=1}^4\frac{2}{3}k_ia_ix_i^{3/2}\sqrt{2g}+\\
%	&&+k_1a_1x_5\sqrt{2gx^*_1}+k_2a_2x_6\sqrt{2gx^*_2}.
%\end{eqnarray}
%where the matrix $J$ is skew-symmetric, and $R$ is positive definite if the controller parameters $k_i$, $i=1\ldots4$,  $k_{I1}$ and $k_{I2}$ are adequately selected (see~\cite{allgoewer} for details).
%
Using this controller the following function:
\begin{small}
\begin{eqnarray}
&H_d(x)&=\frac{1}{2}(x-x^*)^TP^TQP(x-x^*)-u^{*T}Px+\\\notag
&&\sum_{i=1}^4\frac{2}{3}k_ia_ix_i^{3/2}\sqrt{2g}+k_1a_1x_5\sqrt{2gx^*_1}+k_2a_2x_6\sqrt{2gx^*_2},
\end{eqnarray}
\end{small}
which is positive definite and has a global minimum at $x^*$, is an ISS Lyapunov  function with respect to $\varepsilon$, as evidenced by the following bound~\footnote{The expression for the matrix $R$ is not included because of space limitations. The value of $\lambda_m(R)$ can be easily deduced from~\cite{allgoewer}.}:
%$\nabla H_d(x^*)=0$ and 
%$$
%\nabla^2H_d(x)=\textnormal{diag}\lbrace{\frac{k_1a_1\sqrt{2g}}{2\sqrt{x_1}},\frac{k_2a_2\sqrt{2g}}{2\sqrt{x_2}},\frac{k_3a_3\sqrt{2g}}{2\sqrt{x_3}},\frac{k_4a_4\sqrt{2g}}{2\sqrt{x_4}},0,0\rbrace}+P^TQP>0.
%$$ 
%The derivative of $H_d$ along trajectories of the system when measurement errors are present has the form:
%\begin{equation}\notag
%\frac{d}{dt}H_d(\xi)=-\nabla^TH_d(\xi)R\nabla H_d(\xi)-\nabla^TH_d(\xi)g_c'K\varepsilon,
%\end{equation}
%where $g_c'=[g_c^T\;0_{2\times2}]^T$. 
%This shows that if no measurement errors are present $H_d$ is a Lyapunov function for the closed-loop system and the system is asymptotically stable. Furthermore, it also satisfies:
%is an ISS Lyapunov function for $(\xi-x^*)$:
\begin{equation}\notag
\frac{d}{dt}H_d(\xi)\leq-\lambda_m(R)|\nabla H_d(\xi)|^2+|\nabla H_d(\xi)||g_c'K||\varepsilon|.
\end{equation}
This equation suggests the use of the triggering condition:
$$
|\nabla H_d(\xi)||g_c'K||\varepsilon|\leq\rho\lambda_m(R)|\nabla H_d(\xi)|^2, \;\rho\in]0,1[.
$$
%or alternatively:
%$$
%\frac{|\varepsilon|}{|\nabla H_d(\xi)|}\leq(1-\rho)\frac{\lambda_m(R)}{|g_c'K|}, \;\rho\in]0,1[.
%$$
Moreover, assuming the operation of the system to be confined to a compact set containing a neighborhood of $x^*$, $|\nabla H_d(\xi)|$ can be bounded as $|\nabla H_d(\xi)|\geq\rho_m|\xi-x^*|$ and the following triggering rule can be applied to ensure asymptotic stability:
$$
|\varepsilon(t)|^2\leq\sigma|\xi(t)-x^*|^2,\;\sigma=\left(\rho_m\rho\frac{\lambda_m(R)}{|g_c'K|}\right)^2 >0.
$$

\begin{figure}[ht]
%\begin{center}
\centering
\includegraphics[width=\hsize]{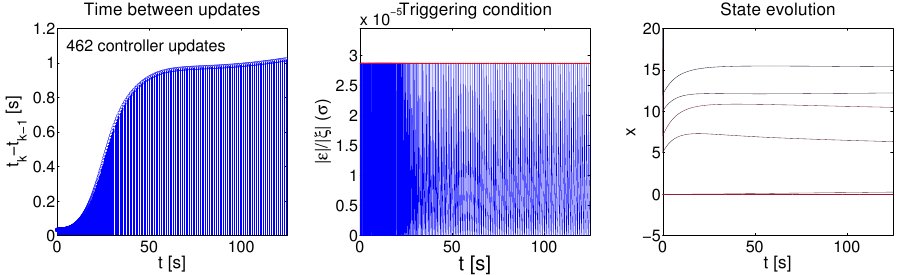}\\%{figure2fix}
\includegraphics[width=\hsize]{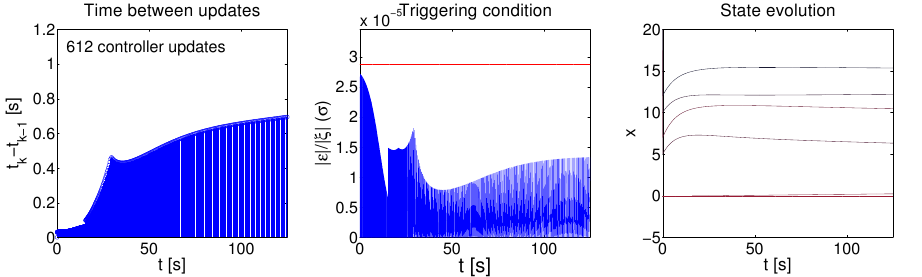}\\%{figure3fix}
\includegraphics[width=\hsize]{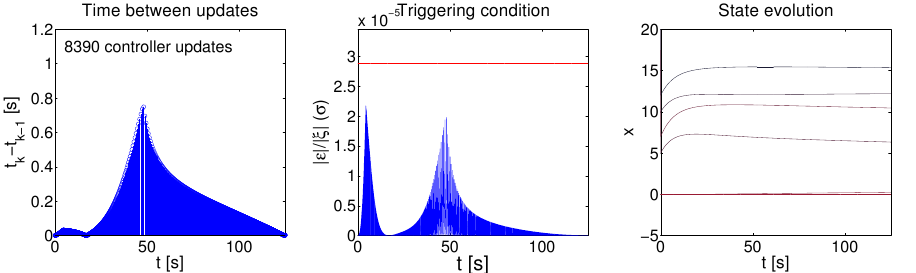}%{figure4fix}
\caption{Times between updates, evolution of the triggering condition, and evolution of the states for the centralized event-triggering implementation (first row), decentralized event-triggering implementation with adaptation (second row), and decentralized event-triggering implementation without adaptation (third row).}
\label{fig:times}
%\end{center}
\end{figure}

\begin{figure}
%\begin{center}
\centering
\includegraphics[width=0.45\hsize]{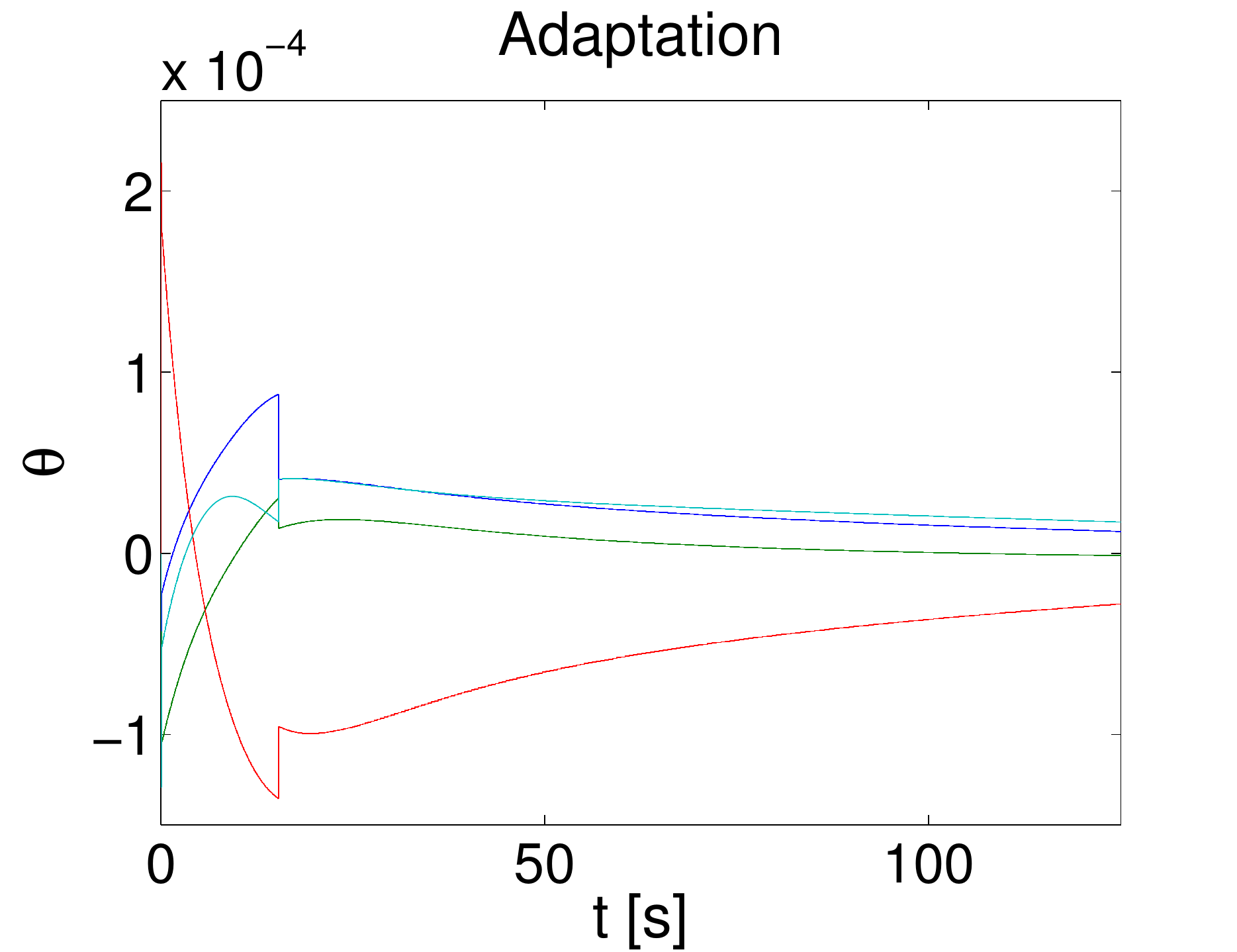}
\caption{Adaptation parameter vector evolution for the adaptive decentralized event-triggered implementation.}
\label{fig:theta}
%\end{center}
\end{figure}

We simulated the decentralized event-triggered implementation of this controller following the techniques in Section~\ref{sec:decentralized_ET}. The physical parameters of the plant and the parameters of the controller are the same as those in~\cite{allgoewer}. Assuming that the system operates in the compact set $S=\lbrace x\in\R^6\;|\;1\leq x_i\leq 20,\; i=1,\dots,\,4;\;0\leq x_i\leq 20,\;i=5,6 \rbrace$, one can take $\rho_m=0.14$, and for the choice of $\rho=0.25$ a value of $\sigma=0.0054^2$ was selected. A bound for the minimum time between controller updates, computed as explained in~\cite{Tabuada}, is given by $\tau_{min}=0.1\,ms$.
%\begin{table}
%%\begin{center}
%\centering
%\begin{tabular}{c c c c c c c }
%  & $A_i$ (cm$^2$) & $a_i$ (cm$^2$) & $k_i$ & $k_{Ii}$ & $Q$ \\
%\hline
%\hline
% $ i=1,2$:	& $50.3$  & $0.233$ & $10$ & $0.13$ & $I_{2\times2}$ \\  
% $ i=3,4$:	& $28.3$  & $0.127$ & $5$ & $-$ & $-$
%\end{tabular}
%%\end{center}
%\caption{Parameter values for the system and controller.\label{tab:params}}
%%\vspace{-1.5cm}
%\end{table}
The decentralized event-triggered controller is implemented adapting $\theta$ as specified by Algorithm~\ref{alg:heuristic} with $q=1$. Furthermore, the pairs of states $x_1,\,x_5$ and $x_2,\,x_6$ are assumed to be measured at the same sensor node, and therefore combined in a single triggering condition at the respective nodes. For comparison purposes, we present in the first row of Figure~\ref{fig:times} the time between controller updates, the evolution of the ratio $\varepsilon/\xi$ \emph{vs} $\sigma$ and the state trajectories, for a centralized event-triggered implementation, starting from initial condition $(12,\,10,\,5,\,7)$ and setting $x_1^*=15$ and $x_2^*=13$. The corresponding results for the proposed decentralized event-triggered implementation are shown in the second row of Figure~\ref{fig:times}, and the results for a decentralized event-triggered implementation without adaptation, \ie, with $\theta(k)=0$ for all $k\in\N$, are shown in the last row of the same figure.
For completeness, Figure~\ref{fig:theta} presents the evolution of adaptation vector $\theta$ for the adaptive decentralized event-triggered implementation.
We can observe that, as expected, a centralized event-triggered implementation is far more efficient, in terms of time between updates, than a decentralized event-triggered implementation without adaption. It is also clear that, although Algorithm~\ref{alg:heuristic} fails to recover the performance of the centralized event-triggered implementation exactly, it produces very good results. The results are even better if we look at the performance in terms of the number of executions which are presented in the legend of these plots. Finally we would like to remark that, although the times between updates in the three implementations can differ quite drastically, the three systems are stabilized producing almost undistinguishable state trajectories.
%The legend in Figures~\ref{fig:times_centralized}, \ref{fig:times_decentralized_adapt} and \ref{fig:times_decentralized_no_adapt} show the number of updates in the figure.  These quantities exemplify the kind of savings obtained by employing the decentralized event-triggered implementation with adaptation compared to a non-adaptive decentralized event-triggered implementation.
%in the $125$ seconds simulated, the adaptive implementation required $59$ controller updates, while the non-adaptive implementation required $622$ controller updates. The centralized event-triggered implementation needed only $43$. 
%*******

\section{Discussion}
\label{sec:discussion}

In~\cite{WangLemmon09} Wang and Lemmon proposed a method for distributed event-triggered control under the assumption that the control system was composed of weakly coupled subsystems. Exploiting this fact, they were able to update inputs independently of each other. Our approach, while not updating inputs independently, does not rely on any internal weak coupling assumptions about the system. Thus, our techniques could be used to complement the techniques in~\cite{WangLemmon09} at the local subsystem level. 

The proposed techniques have been shown effective in decentralizing an event-triggered implementation of a quadruple water-tank system. The centralized controller of this example is a dynamic controller. In this particular case, by allowing the dynamical part of the controller to be continuously computed by the sensors, we successfully obtained a decentralized event-triggered implementation. However, 
the implementation of general dynamic controllers in event-triggered form, centralized or not, remains a question for future research. The design of more efficient adaptation rules is another interesting question to investigate further. Finally, we would like to emphasize the low computational requirements of the proposed implementation, which makes it suitable for sensor/actuator networks with limited computation capabilities at the sensor level.

\bibliographystyle{IEEEtran}
\bibliography{tac_wsn10}

\begin{thebibliography}{10}
\providecommand{\url}[1]{#1}
\csname url@samestyle\endcsname
\providecommand{\newblock}{\relax}
\providecommand{\bibinfo}[2]{#2}
\providecommand{\BIBentrySTDinterwordspacing}{\spaceskip=0pt\relax}
\providecommand{\BIBentryALTinterwordstretchfactor}{4}
\providecommand{\BIBentryALTinterwordspacing}{\spaceskip=\fontdimen2\font plus
\BIBentryALTinterwordstretchfactor\fontdimen3\font minus
  \fontdimen4\font\relax}
\providecommand{\BIBforeignlanguage}[2]{{%
\expandafter\ifx\csname l@#1\endcsname\relax
\typeout{** WARNING: IEEEtran.bst: No hyphenation pattern has been}%
\typeout{** loaded for the language `#1'. Using the pattern for}%
\typeout{** the default language instead.}%
\else
\language=\csname l@#1\endcsname
\fi
#2}}
\providecommand{\BIBdecl}{\relax}
\BIBdecl

\bibitem{wHART_URL}
\BIBentryALTinterwordspacing
{WirelessHART}. [Online]. Available:
  \url{http://www.hartcomm.org/protocol/wihart/wireless_technology.html}
\BIBentrySTDinterwordspacing

\bibitem{RUNES}
K.-E. \AA{}rz\'{e}n, A.~Bicchi, S.~Hailes, K.~Johansson, and J.~Lygeros, ``On
  the design and control of wireless networked embedded systems,'' in
  \emph{2006 IEEE Computer Aided Control System Design,}, Oct. 2006, pp. 440
  --445.

\bibitem{special_proceedings}
P.~Antsaklis and J.~Baillieul, ``Special issue on technology of networked
  control systems,'' \emph{Proceedings of the IEEE}, vol.~95, no.~1, pp. 5 --8,
  Jan. 2007.

\bibitem{Franklin}
G.~Franklin, ``Rational rate [ask the experts],'' \emph{Control Systems
  Magazine, IEEE}, vol.~27, no.~4, pp. 19 --19, Aug. 2007.

\bibitem{Middleton}
G.~Goodwin, R.~Middleton, and H.~Poor, ``High-speed digital signal processing
  and control,'' \emph{Proceedings of the IEEE}, vol.~80, no.~2, pp. 240 --259,
  Feb. 1992.

\bibitem{Goodwin}
G.~Goodwin, S.~Graebe, and M.~Salgado, \emph{Control System Design}.\hskip 1em
  plus 0.5em minus 0.4em\relax Prentice Hall, 2001.

\bibitem{Houpis}
C.~Houpis and G.~B. Lamont, \emph{Digital Control Systems}.\hskip 1em plus
  0.5em minus 0.4em\relax McGraw-Hill Higher Education, 1984.

\bibitem{Nesic:2001p2652}
D.~Nesi{\'c} and A.~Teel, ``Sampled-data control of nonlinear systems: An
  overview of recent results,'' in \emph{Perspectives in robust control}, ser.
  Lecture Notes in Control and Information Sciences, S.~Moheimani, Ed.\hskip
  1em plus 0.5em minus 0.4em\relax Springer Berlin / Heidelberg, 2001, vol.
  268, pp. 221--239.

\bibitem{Nesic:2009}
D.~Nesi{\'c}, A.~Teel, and D.~Carnevale, ``Explicit computation of the sampling
  period in emulation of controllers for nonlinear sampled-data systems,''
  \emph{IEEE Transactions on Automatic Control}, vol.~59, pp. 619--624, Mar.
  2009.

\bibitem{arzen99}
K.-E. {\AA}rz{\'e}n, ``{A simple event based PID controller},'' in
  \emph{Proceedings of 14th IFAC World Congress}, vol.~18, 1999, pp. 423--428.

\bibitem{astrom}
K.~{\AA}str\"{o}m and B.~Bernhardsson, ``Comparison of riemann and lebesgue
  sampling for first order stochastic systems,'' in \emph{Proceedings of the
  41st IEEE Conference on Decision and Control}, vol.~2, Dec. 2002, pp. 2011 --
  2016.

\bibitem{Heemels08}
W.~Heemels, J.~Sandee, and P.~van~den Bosch, ``{Analysis of event-driven
  controllers for linear systems},'' \emph{International Journal of Control},
  vol.~81, no.~4, pp. 571--590, 2008.

\bibitem{Cervin08}
A.~Cervin and T.~Henningsson, ``Scheduling of event-triggered controllers on a
  shared network,'' in \emph{47th IEEE Conference on Decision and Control},
  Dec. 2008, pp. 3601 --3606.

\bibitem{Rabi08}
M.~Rabi and K.~H. Johansson, ``Event-triggered strategies for industrial
  control over wireless networks,'' in \emph{WICON}, 2008.

\bibitem{Rabi082}
M.~Rabi, K.~H. Johansson, and M.~Johansson, ``Optimal stopping for
  event-triggered sensing and actuation,'' in \emph{47th IEEE Conference on
  Decision and Control.}, Dec. 2008, pp. 3607 --3612.

\bibitem{Molin10}
A.~Molin and S.~Hirche, ``Optimal event-triggered control under costly
  observations,'' in \emph{Proceedings of the 19th International Symposium on
  Mathematical Theory of Networks and Systems}, 2010.

\bibitem{Lunze10}
J.~Lunze and D.~Lehmann, ``A state-feedback approach to event-based control,''
  \emph{Automatica}, vol.~46, no.~1, pp. 211 -- 215, 2010.

\bibitem{Tabuada}
P.~Tabuada, ``Event-triggered real-time scheduling of stabilizing control
  tasks,'' \emph{IEEE Transactions on Automatic Control}, vol.~52, no.~9, pp.
  1680 --1685, Sept. 2007.

\bibitem{MazoJr:2008p1340}
M.~Mazo~Jr. and P.~Tabuada, ``On event-triggered and self-triggered control
  over sensor/actuator networks,'' in \emph{Proceedings of the 47th IEEE
  Conference on Decision and Control, 2008.}, Dec. 2008, pp. 435 --440.

\bibitem{Anta:2008p2867}
A.~Anta and P.~Tabuada, ``To sample or not to sample: Self-triggered control
  for nonlinear systems,'' \emph{IEEE Transactions on Automatic Control},
  vol.~55, pp. 2030--2042, Sep. 2010.

\bibitem{MazoAnta}
M.~Mazo~Jr., A.~Anta, and P.~Tabuada, ``An iss self-triggered implementation of
  linear controller,'' \emph{Automatica}, vol.~46, pp. 1310--1314, Aug. 2010.

\bibitem{WangLemmon09}
X.~Wang and M.~D. Lemmon, ``Event-triggering in distributed networked systems
  with data dropouts and delays,'' in \emph{Proceedings of the 12th
  International Conference on Hybrid Systems: Computation and Control}, Apr.
  2009, pp. 366--380.

\bibitem{Soldati}
H.~Zhang, P.~Soldati, and M.~Johansson, ``Optimal link scheduling and channel
  assignment for convergecast in linear wirelesshart networks,'' in
  \emph{Proceedings of the 7th international conference on Modeling and
  Optimization in Mobile, Ad Hoc, and Wireless Networks}, 2009, pp. 82--89.

\bibitem{Sontag:2008p2453}
E.~D. Sontag, ``Input to state stability: Basic concepts and results,'' in
  \emph{Nonlinear and Optimal Control Theory}.\hskip 1em plus 0.5em minus
  0.4em\relax Springer, 2006, pp. 163--220.

\bibitem{allgoewer}
J.~Johnsen and F.~Allg{\"o}wer, ``Interconnection and damping assignment
  passivity-based control of a four-tank system,'' \emph{Lagrangian and
  Hamiltonian Methods for Nonlinear Control 2006}, pp. 111--122, 2007.

\end{thebibliography}

\end{document}